\documentclass[final,1p,times]{elsarticle}

\usepackage{amssymb}
 \usepackage{amsthm}
\usepackage{amscd}
\usepackage{amsmath}
\usepackage{amsfonts}
\usepackage{amssymb}
\usepackage{graphicx}
\newtheorem{theorem}{Theorem}

\newtheorem{corollary}[theorem]{Corollary}

\usepackage{mathrsfs}
\usepackage{titletoc}

\newcommand{\ra}{\rightarrow}
\newcommand{\p}{\partial}
\newcommand{\f}{\frac}

\newcommand{\be}{\begin{equation}}
\renewcommand{\ra}{\rightarrow}
\newcommand{\ee}{\end{equation}}
\newcommand{\bea}{\begin{eqnarray}}
\newcommand{\eea}{\end{eqnarray}}
\newcommand{\bna}{\begin{eqnarray*}}
\newcommand{\ena}{\end{eqnarray*}}

\renewcommand{\le}{\left}
\newcommand{\ri}{\right}

\journal{***}

\begin{document}

\begin{frontmatter}
\title{Global gradient estimate on graph and its applications}

\author{Yong Lin}
 \ead{linyong01@ruc.edu.cn}

\author{Shuang Liu}
 \ead{cherrybu@ruc.edu.cn}

\author{Yunyan Yang}
 \ead{yunyanyang@ruc.edu.cn}
\address{ Department of Mathematics,
Renmin University of China, Beijing 100872, P. R. China}

\begin{abstract}

  Continuing our previous work (arXiv:1509.07981v1), we derive another global gradient estimate for positive functions, particularly for
  positive solutions to the heat equation  on finite or locally finite graphs. In general, the gradient estimate in the present paper is
  independent of our previous one.
  As applications, it can be used to
  get an upper bound and a lower bound of the heat kernel on locally finite graphs.
  These global gradient estimates can be compared with the Li-Yau inequality on graphs contributed by Bauer, Horn, Lin, Lipper, Mangoubi
  and Yau (J. Differential Geom. 99 (2015) 359-409). In many topics, such as eigenvalue estimate and heat kernel estimate
  (not including the Liouville type theorems), replacing
  the Li-Yau inequality  by the global gradient estimate, we can get similar results.
\end{abstract}

\begin{keyword}
 gradient estimate\sep global gradient estimate\sep Harnack inequality\sep locally finite graph

\MSC[2010] 58J35

\end{keyword}

\end{frontmatter}

\titlecontents{section}[0mm]
                       {\vspace{.2\baselineskip}}
                       {\thecontentslabel~\hspace{.5em}}
                        {}
                        {\dotfill\contentspage[{\makebox[0pt][r]{\thecontentspage}}]}
\titlecontents{subsection}[3mm]
                       {\vspace{.2\baselineskip}}
                       {\thecontentslabel~\hspace{.5em}}
                        {}
                       {\dotfill\contentspage[{\makebox[0pt][r]{\thecontentspage}}]}

\setcounter{tocdepth}{2}


\section{Introduction}
The Li-Yau inequality \cite{LiYau2}, or the Li-Yau gradient estimate in literature,  plays an important role in various
topics of geometric analysis, say, the Harnack inequality, the eigenvalue
estimate, the heat kernel estimate, and so on. Recently, this celebrated inequality has been successfully extended from
Riemannian manifolds to finite or locally finite graphs by Bauer, Horn, Lin, Lippner, Mangoubi and Yau \cite{BHLLY}.
Similar to the Riemannian manifold case, the power of the Li-Yau inequality on graphs is evident. For example, it can be
used to get the Liouville type theorems
\cite{BHY}, to estimate eigenvalues and to estimate the heat kernel \cite{BHLLY}.

Note that the Li-Yau inequality is based on the maximum principle, purely local and involving the curvature-dimension condition.
This is exactly the reason why the Li-Yau inequality is so powerful.
Very recently, we find a global gradient estimate for positive functions on finite or
locally finite graphs in \cite{LLY1}. Unlike the Li-Yau inequality, this gradient estimate is based on the graph structure itself,
and thus looks ``rough". As a compensation, it can be applied to positive solutions to
nonlinear differential equations or differential inequalities, such as $\Delta u-\p_tu\leq f(u)$. Due to its global profile,
our estimate can not lead to a Liouville type theorem, since we do not further assume any curvature-dimension condition.
 Indeed, a positive harmonic function is not necessarily constant unless the graph has nonnegative Ricci curvature.
 However, very lucky, we are able to use it to
derive a Harnack inequality, a lower bound of the
first nonzero eigenvalue and an upper bound of the heat kernel \cite{LLY1}. Of course, these bounds are not so delicate as that of \cite{BHLLY},
where the additional curvature-dimension conditions are assumed. \\

In the present paper, exploiting the relation between the Laplacian and the gradient form,
we derive another global gradient estimate, which is independent of our previous one \cite{LLY1}. As applications, we estimate an
upper bound and a lower bound of the heat kernel on locally finite graphs.

\section{Notations and main results}

Suppose that $G=(V,E)$ is a finite or locally finite graph, where $V$ denotes the vertex set and $E$ denotes the edge set.
If $xy\in E$, we assume its weight $w_{xy}>0$. For any $x\in V$, we define the degree of $x$ by
${\rm deg}(x)=\sum_{y\sim x}w_{xy}$. Here and throughout this paper $y\sim x$ if $y$ is adjacent to $x$, or equivalently
$xy\in E$. Let $\mu:V\ra \mathbb{R}$ be a finite measure. Denote $\mu_{\max}=\sup_{x\in V}\mu(x)$,
$w_{\min}=\inf_{x\in V,\,y\sim x}w_{xy}$, $D_\mu=\sup_{x\in V}{{\rm deg}(x)}/{\mu(x)}$, and
$d=\sup_{x\in V,\,y\sim x}{\mu(x)}/{w_{xy}}$. Note that $d\leq \mu_{\max}/w_{\min}$. The
 $\mu$-Laplacian (or the Laplacian for short) acting on a function $f:V\ra \mathbb{R}$ is defined as
$$\Delta f(x)=\f{1}{\mu(x)}\sum_{y\sim x}w_{xy}(f(y)-f(x)).$$
The associated gradient form reads
$$2\Gamma(f,g)(x)=\f{1}{\mu(x)}\sum_{y\sim x}w_{xy}(f(y)-f(x))(g(y)-g(x))$$
for any functions $f:V\ra \mathbb{R}$ and $g:V\ra\mathbb{R}$. Denote $\Gamma(f)=\Gamma(f,f)$.\\

Our main result is the following global gradient estimate

\begin{theorem}\label{Theorem 1}
 Let $G=(V,E)$ be a finite or locally finite graph with $D_\mu<+\infty$. If $u: V\ra\mathbb{R}$ is a positive function, then we
 have
 \be\label{gradient}\f{\Gamma(\sqrt{u})(x)}{u(x)}-\f{\Delta u(x)}{2u(x)}\leq D_\mu,\quad\forall x\in V.\ee
 In particular, if $u:V\times\mathbb{R}\ra \mathbb{R}$ is a positive solution to the heat equation $(\Delta-\p_t)u=0$,
 then we conclude
 \be\label{grad-heat}\f{\Gamma(\sqrt{u})(x)}{u(x)}-\f{\p_t \sqrt{u}(x)}{\sqrt{u}(x)}\leq D_\mu,\quad\forall x\in V.\ee
 \end{theorem}

 In \cite{LLY1}, assuming $D_\mu<+\infty$ and $d<+\infty$, we obtained
 \be\label{pre}\f{\sqrt{2\Gamma(u)}}{u}\leq \sqrt{d}\f{\Delta u}{u}+\sqrt{d}D_\mu+\sqrt{D_\mu}.\ee
 One can check that (\ref{gradient}) is independent of (\ref{pre}).
 Both of them can be easily applied to nonlinear differential inequalities, for examples $\Delta u-\p_tu\leq hu^\alpha$
 and $\Delta u-\p_tu+au\log u\leq 0$,
 where $h$ is a function, $\alpha$ and $a$ are real numbers. The corresponding differential equations
 on Riemannian manifolds were extensively studied,
 see for examples \cite{LiYau1,LiYau2,LiJ,Negrin,Ma,Yang-1,Yang-2} and the references there in.
 In view of Theorem \ref{Theorem 1}, the Li-Yau inequality \cite{BHLLY} can be improved as follows.

 \begin{corollary}\label{cor2} Let $G=(V,E)$ be a finite or locally finite graph with
  $D_\mu<+\infty$ and $D_w=\sup_{x\in V,\,y\sim x}{\rm deg}(x)/w_{xy}<+\infty$.
 Let $u:V\times\mathbb{R}\ra \mathbb{R}$ be a positive function such that $\Delta u-\p_tu=0$ in $B(x_0,2R)$
 for some $R>1$.
 If $G$ satisfies $CDE(n,-K)$ for some $K>0$, then for any $\alpha$, $0<\alpha<1$ and all $t>0$
 $$\f{(1-\alpha)\Gamma(\sqrt{u})}{u}-\f{\p_t\sqrt{u}}{\sqrt{u}}\leq\min\le\{
 D_\mu,\, \f{n}{(1-\alpha)2t}+\f{n(2+D_w)D_\mu}{(1-\alpha)R}+\f{Kn}{2\alpha}
 \ri\}$$
 in the ball $B(x_0,R)$, where $CDE(n,-K)$ is defined as in (\cite{BHLLY}, Definition 3.9; \cite{LinYau}, Definition 1.1).
 \end{corollary}

 As a consequence of Theorem \ref{Theorem 1}, we state the following Harnack inequality.

 \begin{theorem}{\label{Theorem2}} Let $G=(V,E)$ be a finite or locally finite graph with
 $D_\mu<+\infty$, $\mu_{\max}<+\infty$ and $w_{\min}>0$.
 Assume $u:V\times[0,+\infty)\ra \mathbb{R}$ is a positive solution to the heat equation
$\Delta u-\p_tu=0$.
Then for any $(x,T_1)$ and $(y,T_2)$, $T_1<T_2$, we have
\bna
u(x,T_1)\leq u(y,T_2)\exp\le\{2D_\mu(T_2-T_1)+\f{4\mu_{\max}}{w_{\min}}\f{\le({\rm dist}(x,y)\ri)^2}{T_2-T_1}\ri\}.
\ena
\end{theorem}

One of the applications of the Harnack inequality is the estimation of the heat kernel.

 \begin{theorem}{\label{Theorem3}} Let $G=(V,E)$ be a finite or locally finite graph with
 $D_\mu<+\infty$, $\mu_{\max}<+\infty$ and $w_{\min}>0$. Suppose that $w_{xy}=w_{yx}$ for all $x,y\in V$
 with $y\sim x$.
 Let $p(t,x,y)$ be the heat kernel on $G$. Then we state the following:

 \noindent$(i)$ For any $t>0$, $x\in V$ and $y\in V$, there holds
 $$p(t,x,y)\leq \f{1}{{\rm Vol}(B(x,\sqrt{t}))}\exp\le\{
 4\sqrt{\f{2D_\mu \mu_{\max}}{w_{\min}}t}\ri\}.$$
 $(ii)$ If we further assume
 $\mu(x)=\deg(x)$ for all $x\in V$, then
 $$p(t,x,y)\geq \f{1}{\deg(y)}\exp\le\{
-2t-\f{4\mu_{\max}}{w_{\min}}\f{({\rm dist}(x,y))^2}{t}\ri\}.$$
\end{theorem}

For any fixed $x\in V$, it is reasonable that the bound of the heat kernel
implies the volume growth of balls $B(x,\sqrt{t})$ as $t\ra+\infty$. Indeed we have
the following:

\begin{corollary}\label{Cor3'}
Under the assumptions of $(ii)$ of Theorem \ref{Theorem3}, there holds for all $y\in V$ and $t>0$
$${\rm Vol}(B(y,\sqrt{t}))\leq {\rm Vol}(B(y,1))\exp\le\{
t+4\sqrt{\f{2\mu_{\max}}{w_{\min}}t}\ri\}.$$
\end{corollary}

We remark that the above applications have been done in \cite{BHLLY} via the Li-Yau inequality.
By (\cite{LLY1}, Theorem 1), we can get Theorems \ref{Theorem2} and \ref{Theorem3}
with different constants. Also, Theorem \ref{Theorem 1} can lead to (\cite{LLY1}, Theorems 4-6) with different constants.
 The remaining part of this paper
is organized as follows. In Section 3, we prove Theorem \ref{Theorem 1}, Corollary \ref{cor2} and Theorem \ref{Theorem2}.
In Section 4, we prove Theorem \ref{Theorem3} and Corollary \ref{Cor3'}.

\section{Global gradient estimate and Harnack inequality}

In this section, we prove a global gradient estimate (Theorem \ref{Theorem 1}), an improvement of the
Li-Yau inequality (Corollary \ref{cor2}), and a Harnack inequality (Theorem \ref{Theorem2}). We
prove Theorem \ref{Theorem 1} by observing the relation between the gradient form and the Laplacian.
Corollary \ref{cor2} holds just because the terms on the left side of (\ref{grad-heat})
are the same as that of the Li-Yau inequality. For the proof of
Theorem \ref{Theorem2}, we follow the lines of \cite{BHLLY}, and thereby closely follow the lines of \cite{LiYau2}. \\

{\it Proof of Theorem 1}. For any positive function $u: V\ra\mathbb{R}$, there holds
\be\label{delt}2\Gamma(\sqrt{u})=\Delta u-2\sqrt{u}\Delta \sqrt{u}.\ee
 Noting that $\sqrt{u}(y)\geq 0$ for all $y\in V$, and that ${\rm deg}(x)=\sum_{y\sim x}w_{xy}$, we have
\bea\nonumber
-\Delta\sqrt{u}(x)&=&\f{1}{\mu(x)}\sum_{y\sim x}w_{xy}\le(\sqrt{u}(x)-\sqrt{u}(y)\ri)\\\nonumber
&\leq&\f{1}{\mu(x)}\sum_{y\sim x}w_{xy}\sqrt{u}(x)\\\label{g-2}
&=&\f{\deg(x)}{\mu(x)}\sqrt{u}(x).
\eea
Combining (\ref{delt}) and (\ref{g-2}), we conclude (\ref{gradient}).\\

If $u:V\times[0,+\infty)\ra \mathbb{R}$ is a positive solution to the heat equation $\Delta u-\p_t u=0$,
then we have $\Delta u=\p_tu$, and thus (\ref{grad-heat}) follows from (\ref{gradient}) immediately. $\hfill\Box$\\

{\it Proof of Corollary \ref{cor2}}. Let $u:V\ra\mathbb{R}$ be a positive solution to the heat equation $\Delta u=\p_tu$.
By Theorem \ref{Theorem 1},
we have for any $\alpha$, $0<\alpha<1$,
$$\f{(1-\alpha)\Gamma(\sqrt{u})}{u}-\f{\p_t\sqrt{u}}{\sqrt{u}}\leq
\f{\Gamma(\sqrt{u})}{u}-\f{\p_t\sqrt{u}}{\sqrt{u}}\leq D_\mu.$$
This together with the Li-Yau inequality (\cite{BHLLY}, Theorem 4.15) gives the desired result. $\hfill\Box$\\

{\it Proof of Theorem \ref{Theorem2}}. Let $u$ be a positive solution to the heat equation $\Delta u-\p_tu=0$.
By Theorem \ref{Theorem 1}, we have
\be\label{dtlog}-\p_t\log u\leq 2D_\mu-\f{2\Gamma(\sqrt{u})}{u}.\ee
We distinguish two cases to proceed.\\

{\bf Case 1}. $x\sim y$.

For any $s\in [T_1,T_2]$, we have by (\ref{dtlog}) that
\bea\nonumber
\log \f{u(x,T_1)}{u(y,T_2)}
&=&-\int_{T_1}^s\p_t\log u(x,t)dt+\log\f{u(x,s)}{u(y,s)}-
\int_s^{T_2}\p_t\log u(y,t)dt\\
&\leq&2D_\mu(T_2-T_1)-
\int_s^{T_2}\f{2\Gamma(\sqrt{u})(y,t)}{u(y,t)}dt+\log\f{u(x,s)}{u(y,s)}.
\label{temp}\eea
Here we used the fact $-\int_{T_1}^s\f{2\Gamma(\sqrt{u})(x,t)}{u(x,t)}dt\leq 0$.
Note that $x\sim y$, we have
\bea
\nonumber\f{2\Gamma(\sqrt{u})(y,t)}{u(y,t)}&=&\f{1}{\mu(y)}\sum_{z\thicksim y}w_{yz}\le(\f{\sqrt{u}(z,t)}
{\sqrt{u}(y,t)}-1\ri)^2\\
&\geq&\f{w_{\min}}{\mu_{\max}}\le(\f{\sqrt{u}(x,t)}
{\sqrt{u}(y,t)}-1\ri)^2.\label{g1}
\eea
Moreover
\be\label{g2}\log\f{u(x,s)}{u(y,s)}\leq 2\le|\f{\sqrt{u}(x,s)}
{\sqrt{u}(y,s)}-1\ri|.\ee
Inserting (\ref{g1}) and (\ref{g2}) into (\ref{temp}), we obtain
\be\label{inf}
\log \f{u(x,T_1)}{u(y,T_2)}\leq2D_\mu(T_2-T_1)+2\le\{
\le|\f{\sqrt{u}(x,s)}
{\sqrt{u}(y,s)}-1\ri|-\int_{s}^{T_2}\f{w_{\min}}{2\mu_{\max}}\le(\f{\sqrt{u}(x,t)}
{\sqrt{u}(y,t)}-1\ri)^2dt\ri\}.
\ee
Taking the infimum over $s\in [T_1,T_2]$ in (\ref{inf}) and using
(\cite{BHLLY}, Lemma 5.3), we have
$$\log \f{u(x,T_1)}{u(y,T_2)}\leq 2D_\mu(T_2-T_1)+\f{4\mu_{\max}}{w_{\min}}\f{1}{T_2-T_1}.$$

{\bf Case 2}. $x$ is not adjacent to $y$.

Assume ${\rm dist}(x,y)=\ell$. Take a shortest path $x=x_0, x_1,\cdots,x_{\ell}=y$. Let $T_1=t_0<t_1<\cdots
<t_\ell=T_2$, $t_k=t_{k-1}+(T_2-T_1)/\ell$, $k=1,\cdots,\ell$. By the result of Case 1, we have
\bna
\log \f{u(x,T_1)}{u(y,T_2)}&=&\sum_{k=0}^{\ell-1}\le(\log u(x_k,t_k)-\log u(x_{k+1},t_{k+1})\ri)\\
&\leq&\sum_{k=0}^{\ell-1}\le(2D_\mu(t_{k+1}-t_k)+\f{1}{t_{k+1}-t_k}\f{4\mu_{\max}}{w_{\min}}\ri)\\
&\leq& 2D_\mu(T_2-T_1)+\f{\ell^2}{T_2-T_1}\f{4\mu_{\max}}{w_{\min}}\\
&=&2D_\mu(T_2-T_1)+\f{({\rm dist}(x,y))^2}{T_2-T_1}\f{4\mu_{\max}}{w_{\min}}.
\ena
This completes the proof of the theorem.$\hfill\Box$

\section{Heat kernel estimate}

In this section, replacing the Li-Yau inequality in the proof of (\cite{BHLLY}, Theorem 7.6),
we estimate upper bound and lower bound of the heat kernel $p(t,x,y)$. Also, we estimate volume growth
of $B(y,\sqrt{t})$, the ball centered at $y\in V$ with radius $\sqrt{t}$ for any  $t>0$.\\

{\it Proof of Theorem \ref{Theorem3}}. Let $p(t,x,y)$ be the heat kernel on $G$. Let $t>0$ be fixed.
For any $t^\prime>t$ and $x,y,z\in V$, it follows from Theorem \ref{Theorem2} that
$$p(t,x,y)\leq p(t^\prime,z,y)\exp\le\{2D_\mu (t^\prime-t)+\f{4\mu_{\max}}{w_{\min}}
\f{({\rm dist}(x,z))^2}{t^\prime-t}\ri\}.$$
Integrating the above inequality with respect to $z\in B(x,\sqrt{t})$, we obtain
\be\label{ple}p(t,x,y)\leq \f{1}{{\rm Vol}(B(x,\sqrt{t}))}\le(
\sum_{z\in B(x,\sqrt{x})}\mu(z)p(t^\prime,z,y)\ri)\exp\le\{
2D_\mu (t^\prime-t)+\f{4\mu_{\max}}{w_{\min}}
\f{t}{t^\prime-t}\ri\}.\ee
Note that
\be\label{hp1}\sum_{z\in B(x,\sqrt{x})}\mu(z)p(t^\prime,z,y)\leq 1,\ee
and that
\be\label{hp2}\inf_{t^\prime>t}\le(2D_\mu (t^\prime-t)+\f{4\mu_{\max}}{w_{\min}}
\f{t}{t^\prime-t}\ri)=4\sqrt{\f{2D_\mu\mu_{\max}}{w_{\min}}t}.\ee
Inserting (\ref{hp1}) and (\ref{hp2}) into (\ref{ple}), we conclude
$$p(t,x,y)\leq \f{1}{{\rm Vol}(B(x,\sqrt{t}))}\exp\le\{
4\sqrt{\f{2D_\mu\mu_{\max}}{w_{\min}}t}\ri\},$$
and thus $(i)$ holds.

From now on we assume $\mu(x)=\deg(x)$ and thus $D_\mu=1$. Let $p(x,y)=w_{xy}/\deg(x)$, $p_0(x,y)=\delta_{xy}$, where
$\delta_{xy}=1$ if $x=y$, and $\delta_{xy}=0$ if $x\not=y$. Set $p_{k+1}(x,z)=\sum_{y\in V}
p(x,y)p_k(y,z)$. It is well known (see for example \cite{Delmotte}) that the heat kernel $p(t,x,y)$ can be written as
$$p(t,x,y)=e^{-t}\sum_{k=0}^{+\infty}\f{t^k}{k!}\f{p_k(x,y)}{\deg(y)},$$
which leads to
\be\label{abs}p(t,y,y)\geq \f{e^{-t}}{\deg(y)},\quad\forall t>0.\ee
We now distinguish two cases of $t$ to proceed.

{\it Case 1. $t>1$.}

For any $0<\epsilon<1$, it follows from (\ref{abs}) and Theorem \ref{Theorem2} that
$$\f{e^{-\epsilon}}{\deg(y)}\leq p(\epsilon,y,y)\leq p(t,x,y)\exp\le\{
2(t-\epsilon)+\f{4\mu_{\max}}{w_{\min}}\f{({\rm dist}(x,y))^2}{t-\epsilon}\ri\}.$$
Since $t-\epsilon\geq (1-\epsilon)t$, we have
\be\label{low1}
p(t,x,y)\geq \f{e^{-\epsilon}}{\deg(y)}\exp\le\{
-2(t-\epsilon)-\f{4\mu_{\max}}{(1-\epsilon)w_{\min}}\f{({\rm dist}(x,y))^2}{t}\ri\}.
\ee

{\it Case 2. $0<t\leq 1$.}

For any $0<\epsilon<1$, it follows from (\ref{abs}) and Theorem \ref{Theorem2} that
$$\f{e^{-\epsilon t}}{{\deg(y)}}\leq
p(\epsilon t,y,y)\leq p(t,x,y)\exp\le\{2(1-\epsilon)t+\f{4\mu_{\max}}{w_{\min}}\f{({\rm dist}(x,y))^2}{(1-\epsilon)t}\ri\}.$$
And whence in this case
\be\label{low2}p(t,x,y)\geq \f{e^{-\epsilon}}{\deg(y)}\exp\le\{
-2(1-\epsilon)t-\f{4\mu_{\max}}{(1-\epsilon)w_{\min}}\f{({\rm dist}(x,y))^2}{t}\ri\}.
\ee

Combining (\ref{low1}) and (\ref{low2}), we obtain in both cases
$$p(t,x,y)\geq \f{e^{-\epsilon}}{\deg(y)}\exp\le\{
-2t-\f{4\mu_{\max}}{(1-\epsilon)w_{\min}}\f{({\rm dist}(x,y))^2}{t}\ri\}.$$
Letting $\epsilon\ra 0+$, we have
$$p(t,x,y)\geq \f{1}{\deg(y)}\exp\le\{
-2t-\f{4\mu_{\max}}{w_{\min}}\f{({\rm dist}(x,y))^2}{t}\ri\},$$
and thus $(ii)$ holds. $\hfill\Box$\\

{\it Proof of Corollary \ref{Cor3'}}.
By $(i)$ of Theorem \ref{Theorem3}, we have
$$p(t,y,y)\leq \f{1}{{\rm Vol}(B(y,\sqrt{t}))}\exp\le\{
4\sqrt{\f{2D_\mu\mu_{\max}}{w_{\min}}t}\ri\},$$
which together with (\ref{abs}) implies that
$${\rm Vol}(B(y,\sqrt{t}))\leq \deg(y)\exp\le\{
t+4\sqrt{\f{2\mu_{\max}}{w_{\min}}t}\ri\}=
{\rm Vol}(B(y,1))\exp\le\{
t+4\sqrt{\f{2\mu_{\max}}{w_{\min}}t}\ri\}.$$
This gives the desired result. $\hfill\Box$

\bigskip

 {\bf Acknowledgements.} Y. Lin is supported by the National Science Foundation of China (Grant No.11271011). Y. Yang is supported by the National Science Foundation of China (Grant No.11171347 and Grant
 No. 11471014).

\end{document}